\documentclass[fleqn]{mat01}
\usepackage{times,mathtimy,amssymb,latexsym,mathrsfs}%,krepsf,rangecite,mathrsfs}
\begin{document}

\setcounter{page}{531} \firstpage{531}

\def\d{\mbox{\rm d}}
\def\e{\mbox{\rm e}}

%%%%%%%%%%%%%%%%%%%%%%%%%%%%%%
%%%%%%%%%%%  EQUATIONS etc.
%%%%%%%%%%%%%%%%%%%%%%%%%%%%%%
\newcommand{\be}{\begin{equation}}
\newcommand{\ee}{\end{equation}}
\newcommand{\bea}{\begin{eqnarray}}
\newcommand{\eea}{\end{eqnarray}}
\newcommand{\bean}{\begin{eqnarray*}}
\newcommand{\eean}{\end{eqnarray*}}
\newcommand{\brray}{\begin{array}}
\newcommand{\erray}{\end{array}}

\renewcommand{\theequation}{\arabic{section}.\arabic{equation}}

\newtheorem{theore}{Theorem}
\renewcommand\thetheore{\arabic{section}.\arabic{theore}}
\newtheorem{theor}{\bf Theorem}
\newtheorem{lem}[theor]{\it Lemma}
\newtheorem{rem}[theor]{Remark}
\newtheorem{defn}[theor]{\rm DEFINITION}
\newtheorem{exam}[theor]{Example}
\newtheorem{coro}[theor]{\rm COROLLARY}
\newtheorem{propo}[theor]{\rm PROPOSITION}
\newtheorem{remar}{Remark}

%%%%%%%%%%%%%%%%%%%%%%%%%%%%%%%%%
%%%%%%%%%%%% THEOREMS ET AL
%%%%%%%%%%%%%%%%%%%%%%%%%%%%%%%%%
\newtheorem{dfn}{Definition} [section]
\newtheorem{thm}[dfn]{Theorem}
\newtheorem{lmma}[dfn]{Lemma}
\newtheorem{ppsn}[dfn]{Proposition}
\newtheorem{crlre}[dfn]{Corollary}
\newtheorem{xmpl}[dfn]{Example}
\newtheorem{rmrk}[dfn]{Remark}

\newcommand{\bdfn}{\begin{dfn}}
\newcommand{\bthm}{\begin{thm}}
\newcommand{\blmma}{\begin{lmma}}
\newcommand{\bppsn}{\begin{ppsn}}
\newcommand{\bcrlre}{\begin{crlre}}
\newcommand{\bxmpl}{\begin{xmpl}}
\newcommand{\brmrk}{\begin{rmrk}}

\newcommand{\edfn}{\end{dfn}}
\newcommand{\ethm}{\end{thm}}
\newcommand{\elmma}{\end{lmma}}
\newcommand{\eppsn}{\end{ppsn}}
\newcommand{\ecrlre}{\end{crlre}}
\newcommand{\exmpl}{\end{xmpl}}
\newcommand{\ermrk}{\end{rmrk}}

%%%%%%%%%%%%%%%%%%%%%%%%%%%
%%%%%%%%%%%%%%%% SPECIAL SYMBOLS
%%%%%%%%%%%%%%%%%%%%%%%%%%%
\newcommand{\IC}{\mathbb{C}}
\newcommand{\bbc}{\mathbb{C}}
\newcommand{\IZ}{\mathbb{Z}}
\newcommand{\bbz}{\mathbb{Z}}
\newcommand{\IN}{\mathbb{N}}\newcommand {\bbn}{\mathbb{N}}
\newcommand {\bbr}{\mathbb{R}}
%%%%%%%%%%%%%% ABBREVIATIONS
%%%%%%%%%%%%%%%%%%%%%%%%%%%%%%%%

\newcommand{\bta}{\beta}
\newcommand{\eps}{\epsilon}
\newcommand{\LMD}{\Lambda}
\newcommand{\cla}{\mathcal{A}}
\newcommand{\clb}{\mathcal{B}}
\newcommand{\clh}{\mathcal{H}}
\newcommand{\cli}{\mathcal{I}}
\newcommand{\clk}{\mathcal{K}}
\newcommand{\cll}{\mathcal{L}}
\newcommand{\clq}{\mathcal{Q}}
\newcommand{\clu}{\mathcal{U}}
\newcommand{\cls}{\mathcal{S}}
\newcommand{\scrh}{\mathscr{H}}
\newcommand{\scrg}{\mathscr{G}}

\def \bbs {\mbox{\boldmath $s$}}
\def \bbt {\mbox{\boldmath $t$}}

\newcommand{\prf}{\noindent{\it Proof\/}: }
\newcommand{\ots}{\otimes}
\newcommand{\raro}{\rightarrow}
\newcommand{\lraro}{\longrightarrow}
\newcommand{\RARO}{\Rightarrow}
\newcommand{\seq}{\subseteq}
\newcommand{\ol}{\overline}
\newcommand{\lgl}{\langle}
\newcommand{\rgl}{\rangle}
\newcommand{\one}{1\!\!1}
\newcommand{\nn}{\nonumber}

\newcommand{\NI}{\noindent}
\newcommand {\CC}{\centerline}
\def \qed { \mbox{}\hfill
$\Box$\vspace{1ex}}
\newcommand{\half}{\frac{1}{2}}
\newcommand{\id}{\mbox{id}}

\newcommand{\halpha}{\widehat{\alpha}}
\newcommand{\hbeta}{\widehat{\beta}}
%%%%%%%%%%%%%%%%%%%%%%%%%%%%%%%%%
\newcommand{\hta}{\hat{\cla}}
\newcommand{\whtG}{\widehat{G}}
\newcommand{\kernel}{\mbox{ker\,}}
\newcommand{\ran}{\mbox{ran\,}}
\newcommand{\tlds}{\tilde{S}}
%%%%%%%%%%%%%%%%%%%%%%%%%%%%%%%%%
%%%%%%%%%%%%%%%%%%%%%%%%%%%%%%%%%

%%%%%%%%%%%%%%%%%%%%%%%%%%%%%%%%%%%%%%%%

%%\textheight 600pt \textwidth 450pt
%%oddsidemargin 5mm \evensidemargin 5mm \topmargin 0mm

%%\renewcommand{\baselinestretch}{1.2}
%%%\parskip 1ex

\title{On equivariant Dirac operators for $\pmb{SU_q(2)}$}

\markboth{Partha Sarathi Chakraborty and Arupkumar Pal}{On equivariant Dirac operators
for $SU_q(2)$}

\author{PARTHA SARATHI CHAKRABORTY and ARUPKUMAR PAL$^{*}$}

\address{Institute of Mathematical Sciences, CIT Campus, Taramani, Chennai~600~113, India\\
\noindent $^{*}$Indian Statistical Institute, 7, SJSS Marg, New Delhi~110~016, India\\
\noindent E-mail: parthac@imsc.res.in; arup@isid.ac.in\\[1.2pc]
\noindent {\it Dedicated to Prof.~Kalyan Sinha on his
\vspace{-1pc}sixtieth birthday}}

\volume{116}

\mon{November}

\parts{4}

\pubyear{2006}

\Date{}

\begin{abstract}
We explain the notion of minimality for an equivariant spectral triple and show that
the triple for the quantum $SU(2)$ group constructed by Chakraborty and Pal
in~\cite{c-p1} is minimal. We also give a decomposition of the spectral triple
constructed by Dabrowski {\it et~al} \cite{dlssv} in terms of the minimal triple
constructed in~\cite{c-p1}.
\end{abstract}


\maketitle

\section{Introduction}

The interaction between noncommutative geometry and quantum
groups, in particular the (noncommutative) geometry of quantum
groups, had been one of the less understood and less explored
areas of both the theories for a while. In the last few years,
however, there has been some progress in this direction. The first
important step was taken by the authors in \cite{c-p1} where they
found an optimal family of Dirac operators for the quantum $SU(2)$
group acting on $L_2(h)$, the $L_2$ space of the Haar state $h$,
and equivariant with respect to the (co-)action of the group
itself. This family has quite a few remarkable features. They are:

\begin{enumerate}
\renewcommand\labelenumi{\arabic{enumi}.}
\leftskip -.25pc
\item Any element of the $K$-homology group can be realized by a member
from this family, which means that all elements of the $K$-homology group are
realizable through some Dirac operator acting on the single Hilbert space $L_2(h)$ in
a natural manner.

\item The sign of any equivariant Dirac operator on $L_2(h)$ is a compact perturbation of the sign of
a Dirac operator from this family,

\item Given any equivariant Dirac operator $\tilde{D}$
acting on $L_2(h)$, and any Dirac operator $D$ from this family,
there exist
%% one can get a Dirac operator $D$ from this family and%%
two positive reals $k_1$ and $k_2$ such that
\begin{equation*}
|\widetilde{D}| \leq k_1 + k_2|D|.
\end{equation*}
\item
They exhibit features that are unique to the quantum case ($q\neq
1$). It was proved in~\cite{c-p1} that for classical $SU(2)$,
there does not exist any Dirac operator acting on (one copy of)
the $L_2$ space that is both equivariant as well as 3-summable.
\end{enumerate}
These triples were later analysed by Connes \cite{co4} in great
detail, where the general theory of Connes--Moscovici was applied
to obtain a beautiful local index formula for $SU_q(2)$.

Recently, Dabrowski {\it et~al} \cite{dlssv} have constructed
another family of Dirac operators that act on two copies of the
$L_2$ space, has the right summability property, is equivariant in
a sense described in \cite{dlssv}, and is isospectral to the
classical Dirac operator. In this note, we will give a
decomposition of this Dirac operator in terms of the Dirac
operators constructed in \cite{c-p1}.

\section{Equivariance and minimality}

In this section, we will formulate the notion of an equivariant
spectral triple for a compact quantum group and what one means by
its minimality, or irreducibility. For basic notions on compact
quantum groups, we refer the reader to~\cite{w3}. To fix the
notation, let us recall a few things briefly here. Let
$G=(C(G),\Delta)$ be a compact quantum group, where $C(G)$ is the
unital $C^*$-algebra of `continuous functions on $G$' and $\Delta$
the comultiplication map. The symbols $\kappa$ and $h$ will denote
the antipode map and the Haar state for $G$. For two functionals
$\rho$ and $\sigma$ on $C(G)$, the convolution product
$\rho\ast\sigma$ is the functional $a\mapsto
(\rho\otimes\sigma)\Delta(a)$. For $\rho$ as above and $a\in
C(G)$, we will denote by $a\ast\rho$ the element
$(\id\otimes\rho)\Delta(a)$ and by $\rho\ast a$ the element
$(\rho\otimes \id)\Delta(a)$. A unitary representation $u$ of $G$
acting on a Hilbert space $\clh$ is   a unitary element of the
multiplier algebra $M(\clk(\clh)\otimes C(G))$, where $\clk(\clh)$
denotes the space of compact operators on $\clh$, that satisfies
the condition $(\id\otimes\Delta)u=u_{12}u_{13}$. For a unitary
representation $u$ and a continuous linear functional $\rho$ on
$C(G)$, we will denote by $u_\rho$ the operator
$(\id\otimes\rho)u$ on $\clh$. The GNS space associated with the
state $h$ will be denoted by $L_2(h)$ and the cyclic vector will
be denoted by $\Omega$. While using the comultiplication $\Delta$,
we will often use the Sweedler notation (i.e.\
$\Delta(a)=a_{(1)}\otimes a_{(2)}$).

Let $\cla$ be a unital $C^*$-algebra, $G$ be a compact quantum
group and let $\tau$ be an action of $G$ on $\cla$, i.e.\ $\tau$
is a unital $C^*$-homomorphism from $\cla$ to $\cla\otimes C(G)$
satisfying the condition $(\id\otimes\Delta)\tau =
(\tau\otimes\id)\tau$. In other words, let $(\cla, G, \tau)$ be a
$C^*$-dynamical system. Recall \cite{b-s} that a covariant
representation of $(\cla, G,\tau)$ on a Hilbert space $\clh$ is a
pair $(\pi,u)$ where $\pi$ is a unital *-representation of $\cla$
on $\clh$, $u$ is a unitary representation of $G$ on $\clh$ and
they obey the condition
\begin{equation}\label{cov_1}
u(\pi(a)\otimes I)u^* =(\pi\otimes\id)\tau(a),\quad a\in\cla.
\end{equation}
By an {\it odd $G$-equivariant spectral data for $\cla$}, we mean
a quadruple $(\pi,u,\clh,D)$ where
\begin{enumerate}
\renewcommand\labelenumi{\arabic{enumi}.}
\leftskip -.25pc

\item $(\pi,u)$ is a covariant representation of $(\cla, G,\tau)$ on
the Hilbert space $\clh$,

\item $\pi$ is faithful,

\item $u(D\otimes I)u^* =D\otimes I$,

\item $(\pi,\clh,D)$ is an odd spectral triple.
\end{enumerate}
We will often be sloppy and just say $(\pi,\clh,D)$ is an odd
$G$-equivariant spectral triple for $\cla$, omitting $u$. We say
that an operator $D$ on a Hilbert space $\clh$ is an {\it odd
$G$-equivariant Dirac operator} for $\cla$ if there exists a
unitary representation $u$ of $G$ on $\clh$ such that
$(\pi,u,\clh,D)$ gives a $G$-equivariant spectral data for $\cla$.

Similarly, an even $G$-equivariant spectral data for $\cla$
consists of an even spectral data $(\pi,u,\clh, D,\gamma)$ where
$(\pi,u,\clh,D)$ obeys conditions~1, 2 and~3 above, and moreover
$(\pi,\clh,D,\gamma)$ is an even spectral `triple', and one has
$u(\gamma\otimes I)u^* =\gamma \otimes I$. An even $G$-equivariant
Dirac operator is also defined similarly.

We say that an equivariant odd spectral data $(\pi,u,\clh,D)$ is
{\it minimal} if the covariant representation $(\pi,u)$ is
irreducible.

Note that if we take $\cla=C(G)$, then the groups $G$ and $G_{\rm
op}$ have natural actions $\Delta$ and $\Delta_{\rm op}$ on
$\cla$. In what follows, we will mainly be concerned about these
two systems $(\cla=C(G), G, \Delta)$ and $(\cla=C(G), G_{\rm op},
\Delta_{\rm op})$. A $G$-equivariant spectral triple for $C(G)$
will be called a {\it right equivariant spectral triple} for
$C(G)$. A~ {\it right equivariant Dirac operator} for $C(G)$ will
mean a $G$-equivariant Dirac operator for $C(G)$. Similarly, a
$G_{\rm op}$-equivariant spectral triple for $C(G)$ will be called
a {\it left equivariant spectral triple} for $C(G)$ and a $G_{\rm
op}$-equivariant Dirac operator for $C(G)$ will be called a {\it
left equivariant Dirac operator} for $C(G)$.

We will next study covariant representations of the right
$G$-action on $C(G)$, i.e. representations of the system $(C(G),
G,\Delta)$.

%%%%%%%%%%%%%%%%%%%%%%%%%%%%%%%%%%%%%%%%%%%%
\begin{lemma}\label{tech_0}
Let $(\pi, u)$ be a covariant representation of
$(C(G), G, \Delta)$. If the Haar state $h$ of $G$ is faithful,
then $\pi$ is faithful.
\end{lemma}

%%%%%%%%%%%%%%%%%%%%%%%%%%%%%%%%%%%%%%%%%%%%
\begin{proof}
Assume $\pi(a)=0$. Then $\pi(a^*a)=0$ and hence
$(\pi\otimes\id)\Delta(a^*a)=u(\pi(a^*a)\otimes I)u^*=0$. Applying
$(\id\otimes h)$ on both sides, we get $h(a^*a)I=0$. Since $h$ is
faithful,\break $a=0$. \qed
\end{proof}

%%%%%
%% ROLE OF THE ABOVE LEMMA HERE
\begin{remark}{\rm The above lemma helps ensure that if we have a compact
quantum group with a faithful Haar state, take a covariant
representation $(\pi,u)$ of the system $(C(G), G,\Delta)$ on a
Hilbert space $\clh$, and look at a Dirac operator $D$ on $\clh$,
then we really get a spectral triple for the space $G$ rather than
that of some subspace (i.e.\ quotient $C^*$-algebra of
$C(G)$)\break of it.}
\end{remark}

%%%%%%%%%%%%%%%%%%%%%%%%%%%%%%%%%%%%%%%%%%%%
\begin{lemma}\label{tech_1}
Let $(\pi, u)$ be a covariant representation of $(C(G), G,
\Delta)$. Then the operator $u_h$ is a projection and for any
continuous linear functional $\rho$ on $\cla,$ one has $u_h
u_\rho=u_\rho u_h=\rho(1)u_h$.
\end{lemma}
%%%%%%%%%%%%%%%%%%%%%%%%%%%%%%%%%%%%%%%%%%%%

\begin{proof}
Using Peter--Weyl decomposition for $u$, one can assume without
loss in generality that $u$ is finite dimensional. Take two
vectors $w$ and $w'$ in $\clh$. Then
\begin{align*}
\langle w, u_h w'\rangle
 &= \langle w\otimes\Omega , u ( w'\otimes\Omega)\rangle\\[.25pc]
 &= \langle u^*(w\otimes\Omega) ,  w'\otimes\Omega\rangle\\[.25pc]
 &= \langle ((\id\otimes\kappa)u)(w\otimes\Omega) ,
                 w'\otimes\Omega\rangle\\[.25pc]
 &= \overline{\langle w'\otimes\Omega,
      ((\id\otimes\kappa)u)(w\otimes\Omega)\rangle}\\[.25pc]
 &= \overline{\langle w', ((\id\otimes h\kappa)u) w\rangle}\\[.25pc]
 &= \overline{\langle w', ((\id\otimes h)u) w\rangle}\\[.25pc]
 &= \langle u_h w,  w'\rangle.
\end{align*}
Thus $u_h$ is self-adjoint.

Next, for any continuous linear functional $\rho$,
\begin{align*}
u_\rho u_h &= (\id\otimes\rho)u (\id\otimes h)u\\[.25pc]
  &= (\id\otimes\rho\otimes h)(u_{12}u_{13})\\[.25pc]
  &= (\id\otimes\rho\otimes h)(\id\otimes\Delta)u\\[.25pc]
  &= (\id\otimes \rho\ast h)u\\[.25pc]
  &= \rho(1)u_h.
\end{align*}
Similary one has $u_h u_\rho = \rho(1)u_h$.
In particular, $u_h^2=u_h$, so that $u_h$ is a projection.
\qed
\end{proof}

\begin{lemma}\label{tech_2}
Let $A\equiv A(G)$ be the *-subalgebra of $C(G)$ generated by
matrix entries of all finite dimensional unitary representations
of $G$. Let $(A,\clu)$ be a dual pair of Hopf$^{\,\,*}$-algebras
{\rm (}see {\rm \cite{v})}. Then
\begin{equation}\label{dlssv_eq1}
u_\rho \pi(a)=\pi(a\ast \rho_{(1)})u_{\rho_{(2)}}\quad \mbox{for
all }\ \rho\in\clu \mbox{ and } \ a\in A(G).
\end{equation}
\end{lemma}

\begin{proof}
Apply $(\id\otimes\rho)$ on both sides in the equality
\begin{equation*}
u(\pi(a)\otimes I)=((\pi\otimes\id)\Delta(a))u
\end{equation*}
and use the fact that $\rho(ab)=\rho_{(1)}(a)\rho_{(2)}(b)$.\qed
\end{proof}

%%%%%%%%%%%%%%%%%%%%%%%%%%%%%%%%%%%%%%%%%%%%
\begin{lemma}\label{tech_3}
Let $w\in\clh$ be a vector in the range of $u_h$. Then for any
$a\in A(G)$ and $\rho\in\clu${\rm ,} one has
$u_\rho\pi(a)w=\pi(a\ast\rho)w$. In particular{\rm ,} one has $u_h
\pi(a) w= h(a)w$.
\end{lemma}
%%%%%%%%%%%%%%%%%%%%%%%%%%%%%%%%%%%%%%%%%%%%

\begin{proof}
Use Lemma~\ref{tech_2} to get
\begin{align*}
u_\rho\pi(a)w  &= \pi(a\ast \rho_{(1)})u_{\rho_{(2)}}w\\[.2pc]
  &= \pi(a\ast \rho_{(1)})u_{\rho_{(2)}} u_h w\\[.2pc]
  &= \rho_{(2)}(1)\pi(a\ast \rho_{(1)}) u_h w\\[.2pc]
  &= \pi(\rho_{(2)}(1)a\ast \rho_{(1)})  w\\[.2pc]
  &= \pi(a_{(1)}\rho_{(1)}(a_{(2)})\rho_{(2)}(1))w \\[.2pc]
  &= \pi(a_{(1)}\rho(a_{(2)}))w\\[.2pc]
  &= \pi(a\ast \rho)w,
\end{align*}
for $a\in A(G)$.
\qed
\end{proof}

%%%%%%%%%%%%%%%%%%%%%%%%%%%%%%%%%%%%%%%%%%%%
\begin{lemma}\label{tech_4}
The linear span of $\{\pi(a)u_h w\hbox{\rm :}~a\in A(G),
w\in\clh\}$ is dense in $\clh$. In particular{\rm ,} $u_h$ is
nonzero.
\end{lemma}
%%%%%%%%%%%%%%%%%%%%%%%%%%%%%%%%%%%%%%%%%%%%

\begin{proof}
Using Peter--Weyl  decomposition of $u$ and the observation that
$h(\kappa(a))=h(a)$ for all $a\in A$, it follows that
$u_h=(u^*)_h$. Now take a vector $w'$ in $\clh$ such that $\langle
w', \pi(a)u_h w\rangle=0$ for all $w\in\clh$ and $a\in A$. But
then $\langle w', \pi(a)(u^*)_h w\rangle=0$, i.e.\ $\langle
w'\otimes \Omega, (\pi(a)\otimes I)u^* (w\otimes
\Omega)\rangle=0$. The covariance condition~(\ref{cov_1}) now
gives $\langle u (w'\otimes\Omega), (\pi\otimes\id)\Delta(a)
(w\otimes\Omega)\rangle=0$ for all $w\in\clh$ and $a\in A$. In
particular, one has $\langle u (w'\otimes\Omega),
(\pi\otimes\id)\Delta(a) (\pi(b)w\otimes\Omega)\rangle=0$ for all
$w\in\clh$, and $a,b\in A$. Since
$(\pi\otimes\id)\Delta(a)(\pi(b)\otimes
I)=(\pi\otimes\id)(\Delta(a)(b\otimes I))$ and
$\{\Delta(a)(b\otimes I)\hbox{:}~a,b\in A\}$ is total in
$\cla\otimes\cla$, we get $u(w'\otimes\Omega)=0$ and consequently
$w'=0$.\qed
\end{proof}

For $w\in\clh$, denote by $P_w$ the projection onto the closed
linear span of $\{\pi(a)w\hbox{:}~a\in A\}$.

%%%%%%%%%%%%%%%%%%%%%%%%%%%%%%%%%%%%%%%%%%%%
\begin{lemma}\label{tech_5}
Let $w\in u_h\clh$. Then $(P_w\otimes I)u(P_w\otimes
I)=u(P_w\otimes I)$. If $w'$ is another vector in $u_h\clh$ such
that $\langle w,w'\rangle=0${\rm ,} then the projections $P_w$ and
$P_{w'}$ are orthogonal.
% Then $\langle \pi(a)w,\pi(a')w'\rangle=0$
% for any $a,a'\in\cla$.
\end{lemma}
%%%%%%%%%%%%%%%%%%%%%%%%%%%%%%%%%%%%%%%%%%%%

\begin{proof}
For the first part, it is enough to show that $P_w u_\rho P_w
=u_\rho P_w$ for all $\rho\in\clu$. But this is clear because from
Lemma~\ref{tech_3}, we have $u_\rho \pi(a) w = \pi(a\ast \rho) w$.

For the second part, take $a,a'\in A$. Then using
Lemma~\ref{tech_4} one gets
\begin{align*}
\langle \pi(a)w,\pi(a')w'\rangle  &= \langle w,
\pi(a^*a')w'\rangle\\[.2pc]
  &= \langle u_h w, \pi(a^*a')w'\rangle\\[.2pc]
  &= \langle w, u_h \pi(a^*a')w'\rangle\\[.2pc]
  &= \langle w, h(a^*a')w'\rangle\\[.2pc]
  &= 0.
\end{align*}
Thus $P_w$ and $P_{w'}$ are orthogonal.
\qed
\end{proof}

%%%%%%%%%%%%%%%%%%%%%%%%%%%%%%%%%%%%%%%%%%%%
\begin{proposition}$\label{tech_6}\left.\right.$\vspace{.5pc}

\noindent Let $\{w_1,w_2,\ldots\}$ be an orthonormal basis for $u_h\clh$.
Write $P_n$ for $P_{w_n}${\rm ,} and let $\pi_n(\cdot):=P_n
\pi(\cdot)P_n${\rm ,} $u_n := (P_n\otimes I)u(P_n\otimes I)$. Then
\begin{enumerate}
\renewcommand\labelenumi{{\rm \arabic{enumi}.}}
\leftskip -.25pc

\item For each $n${\rm ,} $(\pi_n, u_n)$ is a covariant representation
of the system $(\cla, G,\Delta)$ on $P_n\clh${\rm ,}

\item $\pi=\oplus \pi_n${\rm ,} $u=\oplus u_n${\rm ,}

\item $(\pi_n,u_n)$ is unitarily equivalent to the pair
$(\pi_{\rm L}, u_{\rm R})$ where $\pi_{\rm L}$ is the
representation of $\cla$ on $L_2(G)$ by left multiplications and
$u_{\rm R}$ is the right regular representation of $G$.
\end{enumerate}
\end{proposition}

%%%%%%%%%%%%%%%%%%%%%%%%%%%%%%%%%%%%%%%%%%%%
\begin{proof}
It follows from Lemmas~\ref{tech_5} and \ref{tech_4} that $P_n$'s
are orthogonal, $\sum P_n=I$ and consequently $\pi=\oplus \pi_n$
and $u=\oplus u_n$.

Define $V_n\hbox{:}~P_n\clh\rightarrow L_2(G)$ by
\begin{equation*}
V_n \pi(a)w_n = \pi_L(a)\Omega,\quad a\in A.
\end{equation*}
Since $\langle \pi_{\rm L}(a)\Omega,\pi_{\rm L}(b)\Omega\rangle
=h(a^*b)=\langle \pi(a)w_n,\pi(b) w_n\rangle$,
$\{\pi(a)w_n\hbox{:}~a\in A\}$ is total in $P_n\clh$ and
$\{\pi_{\rm L}(a)\Omega\hbox{:}~a\in A\}$ is total in $L_2(G)$,
$V_n$ extends to a unitary from $P_n\clh$ onto $L_2(G)$. Next, for
$a,b\in A$, one has $V_n\pi(a)\pi(b)w_n=V_n\pi(ab)w_n=\pi_{\rm
L}(ab)\Omega =\pi_{\rm L}(a)\pi_{\rm L}(b)\Omega =\pi_{\rm
L}(a)V_n \pi(b)w_n$. So $V_n \pi(a)=\pi_{\rm L}(a)V_n$ for all
$a\in A$ and hence for all $a\in\cla$.

Finally, we will show that $(V_n\otimes I)u(V_n^*\otimes I) =
u_{\rm R}$. Write $\tilde{u}_n:=(V_n\otimes I)u(V_n^*\otimes I)$.
Then  for any $\rho\in \clu$, one has
\begin{align*}
(\id\otimes \rho)\tilde{u}_n \pi_{\rm L}(a)\Omega
  &= V_n u_\rho V_n^* \pi_{\rm L}(a)\Omega\\[.25pc]
  &= V_n u_\rho \pi(a) w_n \\[.25pc]
  &= V_n u_\rho \pi(a) u_h w_n \\[.25pc]
  &= V_n \pi(a\ast \rho) w_n\\[.25pc]
  &= V_n \pi(a\ast \rho) V_n^* V_n w_n\\[.25pc]
  &= \pi_{\rm L}(a\ast\rho) \Omega.
\end{align*}

By \cite{w3}, $\tilde{u}_n$ must be the right regular
representation $u$ on $L_2(G)$. \qed
\end{proof}

\begin{remark}
{\rm The above proposition leads to an alternative proof of the
Takesaki--Takai duality for compact quantum groups (Theorem~7.5 of
\cite{b-s}).}
\end{remark}
%%%%%%%%%%%%%%%%%%%%%%%%%%%%%%%%%%%%%%%%%%%%
\begin{theorem}[\!]
The covariant representation $(\pi,u)$ is irreducible if and only
if the operator $u_h$ is a rank one projection.
\end{theorem}\pagebreak
%%%%%%%%%%%%%%%%%%%%%%%%%%%%%%%%%%%%%%%%%%%%
\begin{proof}
Immediate corollary of Proposition~\ref{tech_6}. \qed
\end{proof}

\begin{remark}
{\rm In particular, it follows from the above theorem that the
covariant representation $(\pi_{\rm L}, u_{\rm R})$ on $L_2(G)$ is
irreducible. Thus the equivariant Dirac operator constructed
in~\cite{c-p1} is {\it minimal}.}
\end{remark}

%%%%%%%%%%%%%%%%%%%%%%%%%%%%%%%%%%%%%%%%%%%%%%%%%
%%%%%%%%%%%%%%%%%%%%%%%%%%%%%%%%%%%%%%%%%%%%%%%%%
%%%%%%%%%%%%%%%%%%%%%%%%%%%%%%%%%%%%%%%%%%%
%%%%%%%%%%%%%%%%%%%%%%%%%%%%%%%%%%%%%%%%%%%
%%%%%%%%%%%%%%%%%%%%%%%%%%%%%%%%%%%%%%%%%%%

\section{The decomposition}

%%%%%%%%%%%%%%%%%%%%%%%%%%%%%%%%%%%%%%%%%%%
\subsection*{\it Canonical triples for $SU_q(2)$}

Let $q$ be a real number in the interval $(0,1)$. Let $\cla$
denote the $C^*$-algebra of continuous functions on $SU_q(2)$,
which is the universal $C^*$-algebra generated by two elements
$\alpha$ and $\beta$ subject to the relations
\begin{align*}
\alpha^*\alpha+\beta^*\beta&=I=\alpha\alpha^*+q^2\beta\beta^*,
\quad \alpha\beta-q\beta\alpha=0=\alpha\beta^*-q\beta^*\alpha,\\
\beta^*\beta&=\beta\beta^*
\end{align*}
as in \cite{c-p1}. Let $\pi\hbox{:}~\cla\rightarrow\cll(L_2(h))$
be the representation given by left multiplication by elements in
$\cla$. Let $u$ denote the right regular representation of
$SU_q(2)$. Recall~\cite{w3} that $u$ is the unique representation
acting on $L_2(h)$ that obeys the condition
\begin{equation}
((\mbox{id}\otimes\rho)u)\pi(a)\Omega
 = \pi((\mbox{id}\otimes\rho)\Delta(a))\Omega
\end{equation}
for all $a\in\cla$ and for all continuous linear functionals
$\rho$ on $\cla$. In \cite{c-p1}, the authors studied right
equivariant Dirac operators, those Dirac operators that commute
with the right regular representation, i.e.\ $D$ acting on
$L_2(h)$ for which
\begin{equation*}
(D\otimes I)u=u(D\otimes I).
\end{equation*}
In particular, an optimal family of equivariant Dirac operators
were found. A generic member of this family is of the form
\begin{equation*}
e^{(n)}_{ij}\mapsto\begin{cases} (an+b)e^{(n)}_{ij},
&\hbox{if} -n\leq i< n-k,\\[.6pc]%\cr
(cn+d)e^{(n)}_{ij}, &\hbox{if} \ i=n-k,n-k+1,\ldots,n,\end{cases}
\end{equation*}
where $k$ is a fixed nonnegative integer and $a$, $b$, $c$, $d$
are reals with $ac<0$. If one looks at left equivariant Dirac
operators, the same arguments would then lead to the following
theorem.

%%%%%%%%%%%%%%%%%%%%%%%%%%%%%%%
\setcounter{defin}{0}
\begin{theorem}[\!]
Let $v$ be the left regular representation of $SU_q(2)$. Let $k$
be a nonnegative integer and let $a, b, c$, $d$ be real numbers
with $ac<0$. Then the operator $D\equiv D(k,a,b,c,d)$ on $L_2(h)$
given by
\begin{equation*}
e^{(n)}_{ij}\mapsto \begin{cases} (an+b)e^{(n)}_{ij}, &\hbox{if} \
-n\leq j< n-k,\\[.6pc]%\cr
(cn+d)e^{(n)}_{ij}, &\hbox{if}\ j=n-k,n-k+1,\ldots,n,
\end{cases}
\end{equation*}
gives a spectral triple $(\pi,L_2(h),D)$ having nontrivial Chern
character and obeys
\begin{equation}\label{equiv}
(D\otimes I)v=v(D\otimes I).
\end{equation}

Conversely{\rm ,} given any spectral triple $(\pi,L_2(h),
\tilde{D})$ with nontrivial Chern character such that
$(\tilde{D}\otimes I)v=v(\tilde{D}\otimes I)${\rm ,} there exist a
nonnegative integer $k$ and reals $a, b, c$, $d$ with $ac<0$ such
that
\begin{enumerate}
\renewcommand\labelenumi{{\rm \arabic{enumi}.}}
\leftskip -.25pc
\item $\hbox{\rm sign}\,\tilde{D}$ is a compact perturbation of
the sign of $D\equiv D(k,a,b,c,d)$, and

\item there exist constants $k_1$ and $k_2$ such that
\begin{equation*}
|\tilde{D}|\leq k_1+k_2|D|.
\end{equation*}
\end{enumerate}
\end{theorem}

%%%%%%%%%%%%%%%%%%%%%%%%%%%%%%%
\begin{proof}
The key point is to note that the characterizing property of the
left regular representation $v$ is
\begin{equation}
((\mbox{id}\otimes\rho)v^*)\pi(a)\Omega
 = \pi((\rho\otimes\mbox{id})\Delta(a))\Omega.
\end{equation}
Thus on the right-hand side, one now has left convolution of $a$
by $\rho$ instead of right convolution by $\rho$. Therefore any
self-adjoint operator on $L_2(h)$ with discrete spectrum that
obeys $(D\otimes I)v=v(D\otimes I)$ will be of the form
\begin{equation*}
e^{(n)}_{ij}\mapsto \lambda(n,j)e^{(n)}_{ij}.
\end{equation*}
Hence if one now proceeds exactly along the same lines as in
\cite{c-p1}, one gets all the desired conclusions. \qed
\end{proof}

Observe at this point that the whole analysis carried out
in~\cite{co4} will go through for this Dirac operator as well. Let
us now take two such Dirac operators $D_1$ and $D_2$ on $L_2(h)$
given by
\begin{align}
D_1 e^{(n)}_{ij} &=\begin{cases}-2n e^{(n)}_{ij}, &\hbox{if}\
j\neq n\\[.6pc]
(2n+1) e^{(n)}_{ij}, &\hbox{if}\ j=n\end{cases},\nonumber\\[.5pc]
D_2 e^{(n)}_{ij}&=\begin{cases}(-2n-1) e^{(n)}_{ij},
&\hbox{if}\ j\neq n\\[.6pc]
                      (2n+1) e^{(n)}_{ij}, &\hbox{if}\
                      j=n\end{cases}.
\end{align}
Now look at the triple
\begin{equation*}
(\pi\oplus\pi, L_2(h)\oplus L_2(h), D_1\oplus |D_2|).
\end{equation*}
It is easy to see that this is a spectral triple. Nontriviality of
its Chern character is a direct consequence of that of $D_1$. We
will show in the next paragraph that in a certain sense, the
spectral triple constructed in \cite{dlssv} is equivalent to this
above triple.

%%%%%%%%%%%%%%%%%%%%%%%%%%%%%%%%%%%%%%%%%%%%%%%%%%%%
\subsection*{\it The decomposition}

Let us briefly recall the Dirac operator constructed
in~\cite{dlssv}. The carrier Hilbert space $\scrh$ is a direct sum
of two copies of $L_2(h)$ that decomposes as
\begin{equation*}
\scrh=W_0^\uparrow\oplus\left(\bigoplus_{n\in\half\bbz_+}
    ( W_n^\uparrow\oplus W_n^\downarrow )\right),
\end{equation*}
where
\begin{align*}
W_n^\uparrow &=\mbox{span}\bigg\{
u^n_{ij}\hbox{:}~i=-n,-n+1,\ldots,n,\,\\
&\hskip 1.1cm\, j=-n-\half,-n+\half,\ldots,n+\half\bigg\},\\[.4pc]
W_n^\downarrow &=\mbox{span}\bigg\{
d^n_{ij}\hbox{:}~i=-n,-n+1,\ldots,n,\,\\
&\hskip 1.1cm\, j=-n+\half,-n+\frac{3}{2},\ldots,n-\half\bigg\}.
\end{align*}
($u^n_{ij}$ and $d^n_{ij}$ correspond to the basis elements
$|nij\!\uparrow\rangle$ and $|nij\!\downarrow\rangle$ respectively
in the notation of~\cite{dlssv}.) Now write
\begin{equation*}
v^n_{ij}= \begin{pmatrix}u^n_{ij}\\[.5pc]
d^n_{ij}\end{pmatrix}
\end{equation*}
with the convention that $d^n_{ij}=0$ for
$j=\pm\big(n+\half\big)$. Then the  representation $\pi'$ of
$\cla$ on $\scrh$ is given by
\begin{align*}
\pi'(\alpha^{*}) \,v^n_{ij} &=  a^+_{nij} v^{n+\half}_{ i+\half,
j+\half}
 +  a^-_{nij} v^{n-\half}_{i+\half, j+\half},
\nn \\[.5pc]
\pi'(-\beta) \,v^n_{ij} &=  b^+_{nij} v^{n+\half}_{ i+\half,
j-\half}
 +  b^-_{nij} v^{n-\half}_{ i+\half, j-\half},
\nn \\[.5pc]
\pi'(\alpha) \,v^n_{ij} &= \tilde a^+_{nij} v^{n+\half}_{ i-\half,
j-\half}
 + \tilde a^-_{nij} v^{n-\half}_{ i-\half, j-\half},
 \\[.5pc]
 \pi'(-\beta^*) \, v^n_{ij} &= \tilde b^+_{nij} v^{n+\half}_{
i-\half, j+\half}
 + \tilde b^-_{nij} v^{n-\half}_{ i-\half, j+\half},
\nn
\end{align*}
where $a^\pm_{nij}$ and $ b^\pm_{nij}$ are the following  $2
\times 2$ matrices:
\begin{align*}
a^+_{nij} &= q^{\big(i+j-\half\big)/2} [n + i + 1]^\half
\begin{pmatrix}{q^{-n-\half}} \, \frac{[n+j+\frac{3}{2}]^{1/2}}{[2n+2]}
& 0\\[.6pc] q^\half \,\frac{[n-j+\half]^{1/2}}{[2n+1]\,[2n+2]} & q^{-n}
\, \frac{[n+j+\half]^{1/2}}{[2n+1]}\end{pmatrix},
\nn \\[.7pc]
a^-_{nij} &= q^{\big(i+j-\half\big)/2} [n - i]^\half
\begin{pmatrix} {q^{n+1}} \, \frac{[n-j+\half]^{1/2}}{[2n+1]} & -
q^\half
\,\frac{[n+j+\half]^{1/2}}{[2n]\,[2n+1]}\\[.6pc] 0 & q^{n+\half} \,
\frac{[n-j-\half]^{1/2}}{[2n]}\end{pmatrix},
\nn \\[.7pc]
 b^+_{nij} &= q^{\big(i+j-\half\big)/2} [n + i + 1]^\half
\begin{pmatrix}\frac{[n-j+\frac{3}{2}]^{1/2}}{[2n+2]} & 0\\[.6pc] - q^{-n-1}
\,\frac{[n+j+\half]^{1/2}}{[2n+1]\,[2n+2]} & q^{-\half} \,
\frac{[n-j+\half]^{1/2}}{[2n+1]}\end{pmatrix},
\\[.7pc]
b^-_{nij} &=q^{\big(i+j-\half\big)/2} [n - i]^\half\begin{pmatrix}
- q^{-\half} \, \frac{[n+j+\half]^{1/2}}{[2n+1]} & - q^n
\,\frac{[n-j+\half]^{1/2}}{[2n]\,[2n+1]}\\[.6pc] 0 & -
\frac{[n+j-\half]^{1/2}}{[2n]}\end{pmatrix}, \nn
\end{align*}
($[m]$ being the $q$-number $\frac{q^m-q^{-m}}{q-q^{-1}}$) and
$\tilde a^\pm_{nij}$  and $\tilde b^\pm_{nij}$ are the hermitian
conjugates of the above ones:
\begin{equation*}
\tilde a^\pm_{nij} = (a^\mp_{n\pm\half, i-\half, j-\half})^*,
\qquad \tilde b^\pm_{nij} = (b^\mp_{n\pm\half, i-\half,
j+\half})^*.
\end{equation*}
The operator $D$ is given by
\begin{equation*}
D u^n_{ij}=(2n+1) u^n_{ij},
\qquad
D d^n_{ij}= -2n d^n_{ij}.
\end{equation*}
The triple $(\pi', \scrh, D)$ is precisely the triple constructed
in~\cite{dlssv}.

\begin{theorem}[\!]
Let $\clk_q$ be the two-sided ideal of $\cll(\mathscr{H})$
generated by the operator
\begin{equation*}
d^n_{ij}\mapsto q^n d^n_{ij},
\quad
u^n_{ij}\mapsto q^n u^n_{ij},
\end{equation*}
and let $\cla_f$ denote the *-subalgebra of $\cla$ generated by
$\alpha$ and $\beta$. Then there is a unitary $U\hbox{\rm
:}~L_2(h)\oplus L_2(h)\rightarrow \scrh$ such that
\begin{align}
U(D_1\oplus |D_2|)U^{*}  &=   D,\\
U(\pi(a)\oplus\pi(a))U^* -\pi'(a)
     &\in \clk_q \quad \mbox{for all} \ a\in\cla_f.
\end{align}
\end{theorem}

\begin{proof}
Define $U\hbox{:}~L_2(h)\oplus L_2(h)\rightarrow \scrh$ as
follows:
\begin{align*}
U(e^{(n)}_{ij}\oplus 0) &= d^n_{i,j+\half}, \quad
i=-n,-n+1,\ldots,n,\;\\[.4pc]
&\qquad\,j=-n,-n+1,\ldots,n-1,\\[.4pc]
U(e^{(n)}_{in}\oplus 0) &= u^n_{i,n+\half}, \quad
i=-n,-n+1,\ldots,n,\\[.4pc]
U(0\oplus e^{(n)}_{ij}) &= u^n_{i,j-\half}, \quad
i=-n,-n+1,\ldots,n,\; j=-n,-n+1,\ldots,n.
\end{align*}
It is immediate that $U(D_1\oplus |D_2|)U^{*} = D$. Therefore all
that we need to prove now is that $U(\pi(a)\oplus\pi(a))U^{*}
-\pi'(a)\in\clk_q$ for all $a\in \cla_f$. For this, let us
introduce the representation
$\hat{\pi}\hbox{:}~\cla\rightarrow\cll(L_2(h))$ given by
\begin{equation*}
\hat{\pi}(\alpha)=\hat{\alpha}, \qquad
\hat{\pi}(\beta)=\hat{\beta},
\end{equation*}
where $\hat{\alpha}$ and $\hat{\beta}$ are the following operators
on $L_2(h)$ (see Lemma~2.2 of \cite{c-p3}):
\begin{align}
\hat{\alpha}\hbox{:}~e^{(n)}_{ij}  &\mapsto q^{2n+i+j+1}
e^{\big(n+\half\big)}_{i-\half ,j-\half }\nonumber\\[.4pc]
&\quad\,    + (1-q^{2n+2i})^\half(1-q^{2n+2j})^\half
    e^{\big(n-\half\big)}_{i-\half ,j-\half },\label{halpha}\\[.4pc]
\hat{\beta}\hbox{:}~e^{(n)}_{ij}  &\mapsto
   - q^{n+j}(1-q^{2n+2i+2})^\half
            e^{\big(n+\half\big)}_{i+\half ,j-\half }\nonumber\\[.4pc]
&\quad\,  + q^{n+i}(1-q^{2n+2j})^\half
 e^{\big(n-\half\big)}_{i+\half ,j-\half }.\label{hbeta}
\end{align}
It is easy to see that
\begin{equation*}
\pi(a)\oplus\pi(a)-\hat{\pi}(a)\oplus\hat{\pi}(a) \in U^*\clk_q U
\end{equation*}
for $a=\alpha^*$ and $a=\beta$.
Therefore it is enough to verify that
\begin{equation*}
U(\hat{\pi}(a)\oplus\hat{\pi}(a))U^*-\pi'(a)\in\clk_q
\end{equation*}
for $a=\alpha^*$ and for $a=\beta$.

Next observe that
\begin{align*}
 a^+_{nij} &= (1-q^{2n+2i+2})^\half
\begin{pmatrix}
 (1-q^{2n+2j+3})^\half  & 0\\[.5pc]
0 &
 (1-q^{2n+2j+1})^\half
\end{pmatrix} + O(q^{2n}),
 \\[.65pc]
 a^-_{nij} &= q^{2n+i+j+\half}(1-q^{2n-2i})^\half
\begin{pmatrix}
q(1-q^{2n-2j+1})^\half  & 0\\[.5pc] 0 & (1-q^{2n-2j-1})^\half
\end{pmatrix}\\[.5pc]
&\quad\,+O(q^{2n}),
 \\[.65pc]
 b^+_{nij} &= q^{n+j-\half}(1-q^{2n+2i+2})^\half
\begin{pmatrix}
q(1-q^{2n-2j+3})^\half  & 0\\[.5pc] 0 & (1-q^{2n-2j+1})^\half
\end{pmatrix}\\[.5pc]
&\quad\,+O(q^{2n}),
\\[.65pc]
&= q^{n+j-\half}(1-q^{2n+2i+2})^\half
\begin{pmatrix}
q  & 0\\[.2pc] 0 & 1 \end{pmatrix}+O(q^{2n}),
\\[.65pc]
b^-_{nij} &=- q^{n+i}(1-q^{2n-2i})^\half \begin{pmatrix}
(1-q^{2n+2j+1})^\half & 0
\\[.5pc]
0 &(1-q^{2n+2j-1})^\half
\end{pmatrix}\\[.5pc]
&\quad\,+O(q^{2n})\\[.65pc]
&= - q^{n+i} \begin{pmatrix}(1-q^{2n+2j+1})^\half  & 0
\\[.5pc]
0 &(1-q^{2n+2j-1})^\half \end{pmatrix}+O(q^{2n}).
\end{align*}
The required result now follows from this easily.
\qed
\end{proof}

\begin{remark}
{\rm The above decomposition in particular tells us that the
spectral triples $(\pi\oplus\pi, L_2(h)\oplus L_2(h), D_1\oplus
|D_2|)$ and $(\pi', \scrh, D)$ are essentially unitarily
equivalent at the Fredholm module level. Therefore by
Proposition~8.3.14 of \cite{h-r}, they give rise to the same
element in $K$-homology.}
\end{remark}

\begin{remark}
{\rm In the spectral triple in~\cite{dlssv},  the Hilbert space
can be decomposed as a direct sum of two isomorphic copies in such
a manner that in each half Dirac operator has constant sign. So
positive and negative signs come with equal frequency. However
this symmetry is only superficial, as the decomposition above
illustrates. This asymmetry might be a reflection of the inherent
asymmetry in the {\it growth graph} associated with quantum
$SU(2)$ (see \cite{c-p}). For classical $SU(2)$ the graph is
symmetric whereas in the quantum case it is not.

It should also be pointed out here that, at least as far as
classical odd dimensional spaces are concerned, this kind of sign
symmetry is always superficial. They are always inherent  in the
even cases, but not in the odd ones.}
\end{remark}
%%%%%%%%%%%%%%%%%%%%%%%%%%%%%%%%%%%%%%%%%%%%%

\section*{Acknowledgements}
We would like to thank the referee for suggesting some
improvements.

\end{document}